\documentclass[fleqn]{mat01}
\usepackage{times,mathtimy,amssymb}
\begin{document}

\setcounter{page}{33}
\firstpage{33}

\font\zz=msam10 at 10pt
\def\Box{\mbox{\zz{\char'244}}}

\newtheorem{theore}{Theorem}
\renewcommand\thetheore{\arabic{section}.\arabic{theore}}
\newtheorem{theor}[theore]{\bf Theorem}
\newtheorem{propo}[theore]{\rm PROPOSITION}
\newtheorem{lem}[theore]{Lemma}
\newtheorem{definit}[theore]{\rm DEFINITION}
\newtheorem{coro}[theore]{\rm COROLLARY}
\newtheorem{rem}[theore]{Remark}
\newtheorem{exampl}[theore]{Example}
\newtheorem{fact}[theore]{Fact}

\renewcommand{\theequation}{\thesection\arabic{equation}}

\title{Cobordism independence of Grassmann manifolds}

\markboth{Ashish Kumar Das}{Cobordism independence of Grassmann manifolds}

\author{ASHISH KUMAR DAS}

\address{Department of Mathematics, North-Eastern Hill
University, Permanent Campus, Shillong~793~022, India\\
\noindent E-mail: akdas@nehu.ac.in}

\volume{114}

\mon{February}

\parts{1}

\Date{MS received 11 April 2003; revised 9 October 2003}

\begin{abstract}
This note proves that, for $F = \Bbb{R,C}$ or $\Bbb{H}$, the bordism
classes of all non-bounding Grassmannian manifolds $G_k(F^{n+k})$, with
$k < n$ and having real dimension $d$, constitute a linearly independent
set in the unoriented bordism group ${\frak{N}}_d$ regarded as a
${\Bbb{Z}} _2$-vector space.
\end{abstract}


\maketitle

\section{Introduction}\label{S:intro}

This paper is a continuation of the ongoing study of cobordism of
Grassmann manifolds. Let $F$ denote one of the division rings $\Bbb{R}$
of reals, $\Bbb{C}$ of complex numbers, or $\Bbb{H}$ of quaternions. Let
$t = \rm{dim}_{\Bbb{R}}$$F$. Then the Grassmannian manifold
$G_k(F^{n+k})$ is defined to be the set of all $k$-dimensional (left)
subspaces of $F^{n+k}$. $G_k(F^{n+k})$ is a closed manifold of real
dimension $nkt$. Using the orthogonal complement of a subspace one
identifies $G_k(F^{n+k})$ with $G_n(F^{n+k}).$ 

In \cite{pS91}, Sankaran has proved that, for $F = \Bbb{R,C}$ or
$\Bbb{H},$ the Grassmannian manifold $G_k(F^{n+k})$ bounds if and only
if $\nu(n+k) > \nu(k),$ where, given a positive integer $m$, $\nu(m)$
denotes the largest integer such that $2^{\nu(m)}$ divides $m.$ 

Given a positive integer $d$, let $\cal{G}$$(d)$ denote the set of
bordism classes of all non-bounding Grassmannian manifolds
$G_k(F^{n+k})$ having real dimension $d$ such that $k<n.$ The
restriction $k<n$ is imposed because $G_k(F^{n+k}) \approx G_n(F^{n+k})$
and, for $k = n,$ $G_k(F^{n+k})$ bounds. Thus, 
${\cal{G}}(d) = \{ [G_k(F^{n+k})] \in {\frak{N}}_* \mid nkt
= d, k < n,\  {\rm and}\  \nu(n+k) \le \nu(k) \} \subset {\frak{N}}_d.$

The purpose of this paper is to prove the following:

\begin{theor}[\!]\label{T:Gd}
$\cal{G}$$(d)$ is a linearly independent set in the ${\Bbb{Z}}
_2$-vector space ${\frak{N}}_d$.
\end{theor}

Similar results for Dold and Milnor manifolds can be found in
\cite{sK89} and \cite{sD01} respectively.

\section{The real Grassmannians --- a Brief review}\label{S:real}

The real Grassmannian manifold $G_k({\Bbb{R}}^{n+k})$ is an
$nk$-dimensional closed manifold of $k$-planes in ${\Bbb{R}}^{n+k}$. It
is well-known (see \cite{hH80}) that the mod-2 cohomology of
$G_k({\Bbb{R}}^{n+k})$ is given by
\begin{equation*}
H^*(G_k({\Bbb{R}}^{n+k});\ {\Bbb{Z}} _2) \cong {\Bbb{Z}} _2[w_1,
w_2,\dots , w_k, {\bar{w}}_1, {\bar{w}}_2, \dots , {\bar{w}}_n]/
\{w.\bar{w}=1\},
\end{equation*}
where $w = 1 + w_1 + w_2 + \dots + w_k$ and ${\bar{w}} = 1 + {\bar{w}}_1
+ {\bar{w}}_2 +\dots + {\bar{w}}_n$ are the total Stiefel--Whitney
classes of the universal $k$-plane bundle ${\gamma}_k$ and the
corresponding complementary bundle ${\gamma}_k^{\perp},$ both over
$G_k({\Bbb{R}}^{n+k}),$ respectively.

For computational convenience in this cohomology one uses the flag manifold
${\rm Flag}({\Bbb{R}}^{n+k})$  consisting of all ordered $(n+k)$-tuples $(V_1,
V_2,\ldots, V_{n+k})$ of mutually orthogonal one-dimensional subspaces
of ${\Bbb{R}}^{n+k}$ with respect to the `standard' inner product on
${\Bbb{R}}^{n+k}$. It is standard (see \cite{fH66}) that the mod-2
cohomology of ${\rm Flag}({\Bbb{R}}^{n+k})$ is given by
\begin{equation*}
H^*({\rm Flag}({\Bbb{R}}^{n+k});\ {\Bbb{Z}} _2) \cong {\Bbb{Z}} _2[e_1,
e_2, \dots , e_{n+k}] \bigg/ \left\{ \prod_{i=1}^{n+k} (1+e_i) = 1
\right\},
\end{equation*}
where $e_1, e_2,\dots , e_{n+k}$ are one-dimensional classes. In fact
each $e_i$ is the first Stiefel--Whitney class of the line bundle
${\lambda}_i$ over ${\rm Flag}({\Bbb{R}}^{n+k})$ whose total space
consists of pairs, a flag $(V_1, V_2, \ldots, V_{n+k})$ and a vector in
$V_i$.

There is a map ${\pi}_{n+k}\ :\ {\rm Flag}({\Bbb{R}}^{n+k})
\longrightarrow G_k({\Bbb{R}}^{n+k})$ which assigns to $(V_1, V_2, \ldots,
V_{n+k})$, the $k$-dimensional subspace $V_1 \oplus V_2 \oplus \dots
\oplus V_k$. In the cohomology, ${\pi}_{n+k}^*\ :\
H^*(G_k({\Bbb{R}}^{n+k});$ ${\Bbb{Z}}_2) \longrightarrow
H^*({\rm Flag}({\Bbb{R}}^{n+k});\ {\Bbb{Z}}_2)$ is injective and is described
by
\begin{equation*}
{\pi}_{n+k}^*(w) = \underset{i=1}{\overset{k}{\prod}}(1+e_i),\quad  
{\pi}_{n+k}^*(\bar{w}) = \underset{i=k+1}{\overset{n+k}{\prod}}(1+e_i).
\end{equation*}

In \cite{rS82}, Stong has observed, among others, the following facts:

\setcounter{theore}{0}
\begin{fact}\label{F:val}
{\rm The value of the class $u \in H^*(G_k({\Bbb{R}}^{n+k});\ {\Bbb{Z}} _2)$
on the fundamental class of $G_k({\Bbb{R}}^{n+k})$ is the same as the
value of 
\begin{equation*}
{\pi}_{n+k}^*(u)e_1^{k-1}e_2^{k-2} \ldots
e_{k-1}e_{k+1}^{n-1}e_{k+2}^{n-2} \ldots e_{n+k-1}
\end{equation*}
on the fundamental class of ${\rm Flag}({\Bbb{R}}^{n+k})$.}
\end{fact}

\begin{fact}\label{F:for}
{\rm In $H^*({\rm Flag}({\Bbb{R}}^{n+k});\ {\Bbb{Z}}_2)$ one has
\begin{equation*}
e_{i_1}^{n+k-(r-1)}e_{i_2}^{n+k-(r-1)} \dots
e_{i_{r-1}}^{n+k-2}e_{i_r}^{n+k-1}\ =\ 0  
\end{equation*}
if $1 \leq r \leq n+k$ and the set $\{i_1, i_2,\dots , i_r\} \subset
\{1, 2,\dots , n+k\}.$  In particular \ $e_i^{n+k}\ =\ 0$ \ for each
$i$ ,\ $1 \leq i \leq n+k.$}
\end{fact}

\begin{fact}\label{F:top}
{\rm In the top dimensional cohomology of ${\rm Flag}({\Bbb{R}}^{n+k}),$ a monomial
$e_1^{i_1}e_2^{i_2} \dots e_{n+k}^{i_{n+k}}$ represents the non-zero
class if and only if the set $\{i_1 , i_2 , \dots , i_{n+k}\} = \{0 , 1
, \dots , n+k-1\}.$}
\end{fact}

The tangent bundle $\tau$ over $G_k({\Bbb{R}}^{n+k})$ is given (see
\cite{wH64}) by
\begin{equation*}
\tau \oplus {\gamma}_k \otimes {\gamma}_k \cong (n+k){\gamma}_k.
\end{equation*}
In particular, the total Stifel--Whitney class
$W(G_k({\Bbb{R}}^{n+k}))$ of the tangent bundle over $G_k({\Bbb{R}}^{n+k})$
maps under ${\pi}_{n+k}^*$ to
\begin{equation*}
\underset{1 \leq i \leq k}{\prod}
(1+e_i)^{n+k}\  .\   \underset{1
\leq i < j \leq k}{\prod}(1+e_i+e_j)^{-2}.
\end{equation*}
Choosing a positive integer $\alpha$ such that $2^{\alpha} \geq
n+k,$ we have, using Fact \ref{F:for}, 
\begin{equation*}
{\pi}_{n+k}^*(W(G_k({\Bbb{R}}^{n+k}))) = \underset{1 \leq i \leq k}{\prod}
(1+e_i)^{n+k}\  .\   \underset{1
\leq i < j \leq k}{\prod}(1+e_i+e_j)^{2^{\alpha} -2}.
\end{equation*}

$\left.\right.$\vspace{-1.5pc}

\noindent Thus, the $m$th Stiefel--Whitney class $W_m = W_m
(G_k({\Bbb{R}}^{n+k}))$ maps under ${\pi}_{n+k}^*$ to the $m$th elementary
symmetric polynomial in\ $e_i$,\ $1 \leq i \leq k,$\ each with multiplicity
$n+k$,\ and\ $e_i +e_j$,\ $1 \leq i < j \leq k,$\ each with multiplicity\
$2^{\alpha} -2.$  Therefore, if $S_p( {\sigma}_1,{\sigma}_2, \dots ,
{\sigma}_p)$ denotes the expression of the power sum\ $\sum_{m=1}^q y_m^p$ as
a polynomial in elementary symmetric polynomials $\sigma_m$'s in $q$
`unknowns' $y_1,y_2, \dots , y_q$,\ $q \geq p$,\ we have (see \cite{pS91})
\begin{equation*}
S_p ({\pi}_{n+k}^*(W_1),{\pi}_{n+k}^*(W_2), \dots , {\pi}_{n+k}^*(W_p))
= \underset{1 \leq i \leq k}{\sum} (n+k)e_i^p.
\end{equation*}
Thus we have a polynomial
\begin{equation*}
S_p(G_k({\Bbb{R}}^{n+k})) = S_p (W_1,W_2, \dots , W_p) \in
H^p(G_k({\Bbb{R}}^{n+k});\ {\Bbb{Z}} _2)
\end{equation*}
of Stiefel--Whitney classes of $G_k({\Bbb{R}}^{n+k})$ such that 
\setcounter{equation}{3}
\begin{equation}\label{E:sym}
{\pi}_{n+k}^*(S_p(G_k({\Bbb{R}}^{n+k}))) =
\begin{cases}
\underset{1 \leq i \leq k}{\sum} e_i^p,  
&\text{if $n+k$ is odd and\ $p<n+k$}\\
0,  &\text{otherwise.} 
\end{cases} 
\end{equation}

\section{Proof of Theorem~1.1}\label{S:proof}

It is shown in \cite{eF71} that
\begin{equation*}
[G_{2k}({\Bbb{R}}^{2n+2k})] = [G_k({\Bbb{R}}^{n+k})]^4 \quad \text{in}\ 
{\frak{N}}_{4nk}.
\end{equation*}
From this, we have, in particular,
\begin{equation*}
[G_k(F^{n+k})]=[G_k({\Bbb{R}}^{n+k})]^t \quad \text{in}\  {\frak{N}}_{nkt}.
\end{equation*}
For this one has to simply observe that the mod-2 cohomology of the
${\Bbb F}$-Grassmannian is isomorphic as ring to that of the
corresponding real Grassmannian by an obvious isomorphism that
multiplies the degree by $t$. On the other hand, since $\frak{N}_*$ is a
polynomial ring over the field ${\Bbb{Z}} _2$, we have the following:

\setcounter{theore}{0}
\begin{rem}\label{R:ind}{\rm A set
$\{[M_1],[M_2], \dots , [M_m]\}$ is linearly independent in
${\frak{N}}_d$ if and only if the set
$\{[M_1]^{2^\beta},[M_2]^{2^\beta}, \dots , [M_m]^{2^\beta}\}$ is
linearly independent in ${\frak{N}}_{d.2^\beta}$,\ $\beta \geq 0.$}
\end{rem}
Therefore, noting that $t = 1, 2$, or $4$, it is enough to prove
Theorem \ref{T:Gd} for real Grassmannians only. Thus, from now onwards,
we shall take
\begin{equation*}
{\cal{G}}(d) = \{[G_k({\Bbb{R}}^{n+k})]\  \mid \  nk = d,\ k<n,\ \text{and} \;
{\nu}(n+k) \leq {\nu}(k) \}.
\end{equation*}

If $G_k({\Bbb{R}}^{n+k})$ is an odd-dimensional real Grassmannian
manifold then both $n$ and $k$ must be odd, and so ${\nu}(n+k) >
{\nu}(k).$ This means that $G_k({\Bbb{R}}^{n+k})$ bounds and so it
follows that \ ${\cal{G}}(d) = \emptyset$ \ if $d$ is odd. Therefore we
assume that $d$ is even.

\begin{lem}\label{L:Flag}
In $H^*({\rm Flag}({\Bbb{R}}^{n+k});\ {\Bbb{Z}}_2)$ one has{\rm ,} for $1 \leq j
\leq k,$
\begin{align*}
&\left(\underset{1 \leq i \leq k}{\sum} e_i^{n+k-(2j-1)} \right) \cdot
e_1^{k-1}e_2^{k-2} \dots e_{k-j}^j \cdot e_{k-(j-1)}^{j-1} \cdot
e_{k-(j-2)}^{n+k-(j-1)} \dots e_k^{n+k-1} \\[.2pc]
&= e_1^{k-1}e_2^{k-2} \dots e_{k-j}^j \cdot e_{k-(j-1)}^{n+k-j}\cdot
e_{k-(j-2)}^{n+k-(j-1)} \dots e_k^{n+k-1}.
\end{align*}
\end{lem}

\begin{proof}
Note that
\begin{enumerate}
\renewcommand\labelenumi{(\alph{enumi})}
\item if $i \ne k-(j-1)$ then the exponent of $e_i$ in the
product 
\begin{equation*}
\hskip -.5cm e_1^{k-1}e_2^{k-2} \dots
e_{k-j}^j \cdot e_{k-(j-1)}^{j-1} \cdot e_{k-(j-2)}^{n+k-(j-1)} \dots
e_k^{n+k-1}
\end{equation*}
is greater than or equal to $j$, and 

\item $\{n+k-(2j-1)\}+j =n+k-(j-1).$
\end{enumerate}
Therefore, invoking Fact~\ref{F:for}, the lemma follows.\hfill $\Box$
\end{proof}

\begin{propo}\label{P:odd}$\left.\right.$\vspace{.5pc}

\noindent Let ${\cal{O}}(d) = \{[G_k({\Bbb{R}}^{n+k})] \in {\cal{G}}(d)\
\mid \ n+k$ is odd $\}$. Then ${\cal{O}}(d)$ is linearly
independent in~${\frak{N}}_d$.
\end{propo}

\begin{proof}
Arrange the members of ${\cal{O}}(d)$ in descending order of the values
of $n+k$, so that
\begin{equation*}
{\cal{O}}(d) =
\{[G_{k_1}({\Bbb{R}}^{n_1+k_1})],\ [G_{k_2}({\Bbb{R}}^{n_2+k_2})],\ \dots ,
\ [G_{k_s}({\Bbb{R}}^{n_s+k_s})]\},
\end{equation*}
where $n_1+k_1 > n_2+k_2 > \dots > n_s+k_s$.  Note that $n_1 =d$ and
$k_1 = 1$.

For a $d$-dimensional Grassmannian manifold $G_k({\Bbb{R}}^{n+k})$,
consider the polynomials
\begin{equation*}
{f_{\ell}}(G_k({\Bbb{R}}^{n+k})) \; = 
\; \underset{1 \leq j \leq k_{\ell}}{\prod}
S_{n_{\ell} +k_{\ell} -(2j-1)} (G_k({\Bbb{R}}^{n+k}))\; \in
H^d(G_k({\Bbb{R}}^{n+k}); {\Bbb{Z}} _2)
\end{equation*}

$\left.\right.$\vspace{-1.5pc}

\noindent of Stiefel--Whitney classes of $G_k({\Bbb{R}}^{n+k}),$\ where $1
\leq \ell \leq s.$  Then, for each $\ell ,$ $1 \leq \ell \leq s$, we have,
using \eqref{E:sym},
\begin{align*}
&{\pi}_{n_{\ell}+k_{\ell}}^*(f_{\ell}
(G_{k_{\ell}}({\Bbb{R}}^{n_{\ell}+k_{\ell}})))e_1^{k_{\ell}-
1}e_2^{k_{\ell}-2} \dots
e_{k_{\ell}-1}e_{k_{\ell}+1}^{n_{\ell}-1}e_{k_{\ell}+2}^{n_{\ell}-2}
\dots e_{n_{\ell}+k_{\ell}-1} \\[.2pc]
&= \left(\underset{1 \leq j \leq k_{\ell}}{\prod} \left(\underset{1 \leq
i \leq k_{\ell}}{\sum} e_i^{n_{\ell}+k_{\ell}-(2j-1)}
\right)\right)e_1^{k_{\ell}-1}e_2^{k_{\ell}-2} \dots
e_{k_{\ell}-1}e_{k_{\ell}+1}^{n_{\ell}-1}e_{k_{\ell}+2}^{n_{\ell}-2}\\[.2pc]
&\quad\ \dots e_{n_{\ell}+k_{\ell}-1} \\[.2pc]
&= e_1^{n_{\ell}}e_2^{n_{\ell}+1}
\dots e_{k_{\ell}}^{n_{\ell}+k_{\ell}-1}e_{k_{\ell}+1}^{n_{\ell}-1}
e_{k_{\ell}+2}^{n_{\ell}-2} \dots e_{n_{\ell}+k_{\ell}-1}, 
\end{align*}
applying Lemma \ref{L:Flag} repeatedly for successive values of $j$.

Thus, in view of Facts~\ref{F:val} and \ref{F:top}, the Stiefel--Whitney
number
\begin{equation*}
\langle {f_{\ell}}(G_{k_{\ell}}({\Bbb{R}}^{n_{\ell}+k_{\ell}})),
[G_{k_{\ell}}({\Bbb{R}}^{n_{\ell}+k_{\ell}})] \rangle \; \ne \; 0
\end{equation*}
for each $\ell$, $1 \le \ell \le s.$  On the other hand, using
\eqref{E:sym},  it is clear that
\begin{equation*}
\langle {f_{\ell}}(G_{k_h}({\Bbb{R}}^{n_h+k_h})),
[G_{k_h}({\Bbb{R}}^{n_h+k_h})] \rangle \; = \; 0
\end{equation*}
for each $h > \ell$, since $n_{\ell}+k_{\ell}-1 \; \ge \; n_h+k_h.$
Therefore, it follows that the $s \times s$ matrix
\begin{equation*}
[\langle {f_{\ell}}(G_{k_h}({\Bbb{R}}^{n_h+k_h})),
[G_{k_h}({\Bbb{R}}^{n_h+k_h})] \rangle ]_{1 \le {\ell}\le s,\; 1 \le h
\le s}
\end{equation*}
is non-singular; being lower triangular with 1's in the diagonal.  This
completes the\break proof.\hfill $\Box$
\end{proof}

Now we shall complete the proof of Theorem \ref{T:Gd} using induction on $d$. 
First note that
\begin{align*}
{\cal G}(2) &= \{[G_1({\Bbb{R}}^{2+1})]\} = \{[{\Bbb{R}}P^2]\},\\[.2pc]
{\cal G}(4) &= \{[G_1({\Bbb{R}}^{4+1})]\} = \{[{\Bbb{R}}P^4]\},
\end{align*}
and so both are linearly independent in ${\frak{N}}_2,$ ${\frak{N}}_4$
respectively. Assume that the theorem holds for all dimensions less than
$d$.

We have ${\cal G}(d) \; = \; {\cal E}(d) \bigcup {\cal O}(d) $,\; where
\begin{equation*}
{\cal E}(d) \; = \; \{[G_k({\Bbb{R}}^{n+k})] \in {\cal G}(d) | n+k\ \ 
\text{is even}\}
\end{equation*}
and
\begin{equation*}
{\cal O}(d) \; = \; \{[G_k({\Bbb{R}}^{n+k})] \in
{\cal G}(d) | n+k\ \ \text{is odd}\}.
\end{equation*}
Observe that if $[G_k({\Bbb{R}}^{n+k})] \in {\cal E}(d)$ then both $n$
and $k$ are even with $\nu(k) \ne \nu(n).$ On the other hand,
$[G_2({\Bbb{R}}^{\frac{d}{2} +2})] \ \in \ {\cal E}(d)$ \ if $ d \equiv
0 \pmod{8}.$ Thus, ${\cal E}(d) \ne \emptyset $ if and only if $ d
\equiv 0 \pmod{8}.$

In view of Proposition \ref{P:odd}, we may assume without any loss that
${\cal E}(d) \; \ne \; \emptyset .$ Then, by the above observation and
by Theorem 2.2 of \cite{pS91} every member of ${\cal E}(d)$ is of the
form $[G_{\frac{k}{2}}({\Bbb{R}}^{\frac{n}{2} +\frac{k}{2}})]^4$,
where $[G_{\frac{k}{2}}({\Bbb{R}}^{\frac{n}{2} +\frac{k}{2}})]$ $\in$ $
{\cal G}(\frac{d}{4}).$ By induction hypothesis, $ {\cal
G}(\frac{d}{4})$ is linearly independent in ${\frak{N}}_{\frac{d}{4}}$. 

So, by Remark~\ref{R:ind},
\setcounter{equation}{3}
\begin{equation}\label{E:ind}
{\cal E}(d) \ {\rm is\ linearly\ independent\ in} \ {\frak{N}}_d.
\end{equation}

Again note that if $[G_k({\Bbb{R}}^{n+k})] \; \in \; {\cal E}(d),$ then,
by (\ref{E:sym}), the polynomial $S_p(G_k({\Bbb{R}}^{n+k}))$ $=$ $0, \;
\forall \; p \ge 1$. So, for each of the polynomials $f_{\ell},$\; $1
\le \ell \le s,$ considered in Proposition \ref{P:odd}, we have
\begin{equation*}
\langle {f_{\ell}}(G_{k}({\Bbb{R}}^{n+k})), [G_{k}({\Bbb{R}}^{n+k})]
\rangle \; = \; 0.
\end{equation*}
Therefore, writing 
\begin{align*}
\hskip -4pc {\cal E}(d) = \{ [G_{k_{s+1}}({\Bbb{R}}^{n_{s+1}+k_{s+1}})],
[G_{k_{s+2}}({\Bbb{R}}^{n_{s+2}+k_{s+2}})],\dots,
[G_{k_{s+q}}({\Bbb{R}}^{n_{s+q}+k_{s+q}})]\},
\end{align*}
where $n_{s+1}+k_{s+1} > n_{s+2}+k_{s+2} > \dots > n_{s+q}+k_{s+q},$ we
see that the $s \times (s+q)$ matrix 
\begin{equation*}
[\langle {f_{\ell}}(G_{k_h}({\Bbb{R}}^{n_h+k_h})),
[G_{k_h}({\Bbb{R}}^{n_h+k_h})] \rangle ]_{1 \le {\ell}\le s,\; 1 \le h
\le s+q}
\end{equation*}
is of the form
\begin{equation}\label{M:od}
\underset{{\cal O}(d) \qquad \qquad \qquad \qquad {\cal E}(d)}{
\left[\begin{matrix}
1 & 0 & 0 & \cdots & 0 \quad & 0 & 0 & 0 &\cdots & 0\\
 \star & 1 & 0 & \cdots & 0 \quad & 0 & 0 & 0 &\cdots & 0\\
\star & \star & 1 & \cdots & 0 \quad & 0 & 0 & 0 &\cdots & 0\\
 - & - & - & \cdots & - \quad & - & - & - &\cdots & -\\
 - & - & - & \cdots & - \quad & - & - & - &\cdots & -\\
\star & \star & \star & \cdots & 1 \quad & 0 & 0 & 0 &\cdots & 0
\end{matrix}
\right]}.
\end{equation}
Thus, no non-trivial linear combination of members of ${\cal O}(d)$ can
be expressed as a linear combination of the members of ${\cal E}(d).$
This, together with (\ref{E:ind}) and Proposition~\ref{P:odd}, proves
that the set ${\cal G}(d) \; = \; {\cal E}(d) \bigcup {\cal O}(d) $\ is
linearly independent in ${\frak{N}}_d$. Hence, by induction, 
Theorem~\ref{T:Gd} is completely proved.

\setcounter{theore}{5}
\begin{rem}{\rm
Using the decomposition of the members of ${\cal E}(d)$, and the
polynomials $f_{\ell}$, in the lower dimensions together with the {\it
doubling homomorphism} defined by Milnor \cite{jM65}, one can obtain a
set of polynomials of Stiefel--Whitney classes which yield, as in
Proposition~\ref{P:odd}, a lower triangular matrix for ${\cal E}(d)$
with $1$'s in the diagonal. Thus using (\ref{M:od}) we have a lower
triangular matrix, with 1's in the diagonal, for the whole set ${\cal
G}(d).$}
\end{rem}

\section*{Acknowledgement}

Part of this work was done under a DST project

\end{document}